\documentclass[11pt, reqno, a4paper]{amsart}
\usepackage[margin=1.3in]{geometry}
\numberwithin{equation}{section}
\usepackage{amssymb,amsfonts,amsthm}
\usepackage[utf8]{inputenc}
\usepackage{listings}
\usepackage{bm}
\usepackage[hyperpageref]{backref}

\usepackage{color}
\usepackage[bookmarks]{hyperref}

\newcommand{\dd}{\, \mathrm{d}}
\newcommand{\R}{\ensuremath{\mathbb R}}

\allowdisplaybreaks[4]

\newtheoremstyle{myremark}{10pt}{10pt}{}{}{\bfseries}{.}{.5em}{}

 \newtheorem{thm}{Theorem}[section]
 \newtheorem{cor}[thm]{Corollary}
 \newtheorem{lem}[thm]{Lemma}
 \newtheorem{prop}[thm]{Proposition}
 \theoremstyle{definition}

 \newtheorem{rem}[thm]{Remark}

\newcommand{\Om} {\Omega}

 \newcommand{\C}{\mathbb{C}}

 \newcommand{\norm}[1]{\left\Vert#1\right\Vert}

\def\Dp{{{\mathcal D}^{1,p}_0(\Om)}}
\def\Dt{{{\mathcal D}^{1,2}_0(\Om)}}
\def\w2r{{{ W}^{2,2}(\R^N)}}

\def\cset{{\subset \subset }}

\def\C{{\mathcal C}}
\def\D{{\mathcal D}}

\def\H{{\mathcal{H}_p(\Om)}}

\def\R{{\mathbb R}}

\def\F{{\mathcal F}}
\def\({{\Big(}}
\def\){{\Big)}}

\def\ws2{{\F_{\frac{N}{2}}}}
\def\L2{{ L^{1,\;\infty}(\log L)^2}}
\def\l2{\mathcal M\log L}

\def\c1Loc{{\C_{loc}^1}}
\def\cp{{\rm Cap}_p}

\allowdisplaybreaks[4]

\begin{document}

\title[Hardy inequality ]{Quantitative stability for weighted Hardy-type inequalities}


 \author[S. Roy]{ Subhajit Roy}
\address[Subhajit Roy]{Department of Mathematics, Indian Institute of Science Education and Research, Pune 411008, India.}
\email{ rsubhajit.math@gmail.com, subhajit.roy@iiserpune.ac.in}

\thanks{}

\maketitle

\begin{abstract} For $p\in(1,\infty)$ and a domain $\Omega\subseteq\mathbb R^N$ (bounded when $p\ge N$), we consider the following weighted Hardy-type inequality:
\begin{equation}\label{asb}
\int_\Omega|\nabla u|^p\dd x\ge C_H \int_\Omega|g||u|^p\dd x,\quad u\in \Dp,   
\end{equation}
 where $C_H$ is the optimal constant in \eqref{asb}, and $\Dp$ denotes the Beppo-Levi  space. In this article, we investigate {\it{quantitative stability}} (including {\it{quantitative gradient stability}}) for \eqref{asb} under suitable assumptions on the weight $g\in L_{\mathrm{loc}}^1(\Omega)$. More precisely, we identify a class of weights $g$ for which quantitative stability estimates hold. In the linear case $p=2$, we also study the {\it{optimal stability constants}} in the gradient stability estimates and whether they are {\it{attained}}. Furthermore, we prove that {\it{the exponent}} in the gradient stability estimate is optimal for certain weights $g$.

 \end{abstract}

\section{Introduction}
There has been a growing interest in studying the stability for various functional and geometric inequalities, which play a crucial role in several problems arising in the calculus of variations, partial differential equations, and geometry, etc. 
Let
$\mathcal{E}(u)\geq \mathcal{F}(u)$
be a functional inequality on a function space $X$, where $\mathcal{E}(\cdot), \mathcal{F}(\cdot)$ are certain energy functionals. 
Assume that the set of extremal functions 
$\mathcal{M}:=\{u\in X:\mathcal{E}(u)=\mathcal{F}(u)\}$ is non-empty. A fundamental question in the theory of functional inequalities is whether or not the inequality admits a \emph{quantitative stability estimate}. More precisely, the question is: if $\mathcal{E}(u)-\mathcal{F}(u)$ is very close to an equality, i.e., the {\it deficit functional}  $\delta(u):=\mathcal{E}(u)-\mathcal{F}(u)$ is close to $0$, then is it necessary that 
$u$ has to be `close' to the set of extremals $\mathcal{M}$? Such a property is typically expressed by a quantitative estimate of the form
\begin{align} \label{Eq:gen_stab}
    \delta(u)\geq C_S\,d(u,\mathcal{M}) \ \ \ \forall u \in X,
\end{align}
where $C_S>0$ is a constant and $d(u,\mathcal{M})$ is a suitable notion of distance from $u$ to $\mathcal{M}$.
The question of quantitative stability for the $L^2$-Sobolev inequality was first raised by Brezis and Lieb \cite{Brezis-Lieb} and was settled by Bianchi and Egnell \cite{Bianchi}.
They showed that
there exists $C_S>0$ such that for every $u\in W^{1,2}_0(\mathbb{R}^N),$
\begin{align} \label{Eq:Sob_stab}
     \|\nabla u\|_{L^2(\mathbb{R}^N)}^2-\mathcal{S}_N^2\|u\|_{L^{2^*}(\mathbb{R}^N)}^2\geq C_S\,\inf_{v\in\mathcal{M}}\|\nabla u-\nabla v\|_{L^2(\mathbb{R}^N)}^2,
\end{align} where $N\ge 3,\,2^*=\frac{2N}{N-2},$ $\mathcal{S}_N$ is the optimal Sobolev constant, and $\|\cdot\|_{L^2(\R^N)}$ denotes the $L^2$-norm. Since then, quantitative stability for the $L^p$-Sobolev inequality \cite{Ci-Fus2009, Fig-New2019,Neu2020, Fi-Zh2022, Bhakta}, and more generally functional inequalities, has become an active area of research.  
The literature in this direction is vast, and the references below are by no means exhaustive. We mention a few works related to our objectives, e.g., on Poincar\'e inequality \cite{Lu-poin,Amir2},  geometric Hardy inequality \cite{Roy-Das}, Gagliardo-Nirenberg-Sobolev inequality \cite{carlen}, Caffarelli-Kohn-Nirenberg inequality \cite{wei-caf, Anh-Lu}, and log-Sobolev inequality \cite{Jean,Chen}.  
{ Quantitative stability of the classical Hardy inequality with respect to the distance from its {\it virtual extremals} was established in \cite{Chianchi_2008}. 

\smallskip

A quantitative estimate like \eqref{Eq:gen_stab} naturally raises the question of determining the optimal constant $C_s$ and whether it is attained. 
In the linear case $p=2,$ the attainment and lower bound of {\it{best stability constant}} is studied for the geometric Hardy inequality \cite{Roy-Das},  Sobolev inequality \cite{ konig, Jean}, Hardy-Sobolev inequality \cite{Souptik}, and Caffarelli-Kohn-Nirenberg inequality \cite{Cazacu, wei-Wu}; see the references therein for further development on this topic. 

\smallskip

Let 
$\Omega$ be a domain in $\R^N$ such that $\Omega\;\text{{\it{is bounded}} if}\;p\ge N,$ while $\Omega$ can be any domain if $p\in (1,N).$ 
Consider the {\it{``Beppo-Levi" space}} $\Dp$ defined as the closure of $C_c^{\infty}(\Omega)$ with respect to the norm, $ \norm{u}_{\Dp} :=\|\nabla u\|_{L^p(\Omega)}.$  
Clearly, $\Dt$ is an inner product space with inner product $\langle u,v\rangle_{\Dt}:=\int_\Omega\nabla u\cdot\nabla v\dd x$ for $u,v\in \Dt.$
 The aim of this paper is to investigate  quantitative (gradient) stability of the following weighted Hardy-type inequality
 \begin{eqnarray}\label{HS}
 \|\nabla u\|_{L^p(\Omega)} \ge C_H \| u\|_{L^p(\Omega,\,|g|)}  , \quad\forall\; u \in \Dp,
 \end{eqnarray}
for some weight function $g\in L^1_{\mathrm{loc}}(\Omega)$ for which $C_H$ is attained in $\Dp$, where $C_H>0$ is the optimal constant in \eqref{HS} and 
$$\|u\|_{L^p(\Omega,\,|g|)}:=\left(\int_{\Omega} |g||u|^p\dd x\right)^{1/p}.$$ 
We call a weight function $g\in L^1_{\mathrm{loc}}(\Omega)$ Hardy potentials if $g$ satisfies \eqref{HS}.
We denote the space of the Hardy potentials by $\H.$

\smallskip

 Using Poincar\'{e} inequality, it is easy to verify that $L^{\infty}(\Omega)\subseteq \H$ if $\Omega$ is bounded (in one direction). Moreover, for a bounded convex domain $\Omega$, the geometric Hardy inequality
  \begin{equation*}
 \|\nabla u\|_{L^p(\Omega)} \ge (p-1)/p\; \| u\|_{L^p(\Omega,\,d_\Omega^{-p})}   \quad \forall\,u \in \Dp,
  \end{equation*}
  ensures that $d_\Omega^{-p}\in \H,$ where $d_\Omega(x)=\inf\{|x-y|:y\in \partial\Omega\}.$ Several other function spaces are also known to be contained in $\H$. For $p=2$ and bounded $\Omega$, it is known that $L^{N/2}(\Omega)\subseteq \H$ \cite{alleg}. For $p\in(1,N)$, the Lorentz–Sobolev embedding gives $L^{N/p,\infty}(\Omega)\subseteq \H$ \cite{Vis-05}. Furthermore, for bounded domains $\Omega$, it was shown in \cite{Nirjan} that $L^q(\Omega)\subseteq\H$ for $q>1$ when $p=N$, and for $q=1$ when $p>N$. See also \cite{subha1,Roy-a} and the references therein for a more general class of admissible function spaces contained in $\H$.

  \smallskip
  
   In \cite{Mazya}, Maz'ya gave a very intrinsic characterization of a Hardy potential using the $p$-capacity (see Section 2.4.1, page 128).
 Recall that, for $ F \cset \Omega,$ the $p$-capacity of $F$ relative to $\Omega$ is defined as,
  \[{\cp(F,\Omega)}  = \inf \left\{ \displaystyle \int_{\Omega} | \nabla u |^p \ {\mathrm{d}}x : u \in \Dp : u \geq 1 \;\text{on}\;F \right \}.\]
  Maz'ya proved that  $g\in\H$ if and only if
  \begin{eqnarray*}
  \norm{g}_{\H}:= \sup\left\{ \frac{ \int_{F} |g| \ {\mathrm{d}}x}{\cp(F,\Omega)}:F \cset \Omega; |F|\ne 0 \right\}<\infty.                           
   \end{eqnarray*}  
One can verify that $\H=\{g\in \H: \|g\|_{\H}<\infty\}$ and $\norm{.}_\H$ is a Banach function norm on $\H $. 

 \smallskip
 
  For $g\in \H,$ let $\mathbb{H}_{p}(\Omega)$ be the best constant  in \eqref{HS} (Defined in  \eqref{best-const}). 
One of the sufficient conditions ensuring that $\mathbb{H}_{p}(\Omega)$ is attained in $\Dp\setminus\{0\}$ (i.e., equality holds in \eqref{HS}) is the compactness of the map
 $$G_p(u):= \displaystyle \int_{\Om} |g||u|^p \ {\mathrm{d}}x\,,\quad u\in \Dp$$   
(i.e.,  $u_n\rightharpoonup u$ implies $G_p(u_n)\to G_p(u)$). This can be proved using variational methods.
Indeed, if $G_p$ is compact, then $u$ is an extremal for \eqref{HS} if and only if $u$ belongs to the family 
$$ \mathcal{M}_p:=\{c\,\Phi: c\in \mathbb{R}\}$$
where  $\Phi\in \Dp\setminus\{0\}$ is a non-negative extremal for \eqref{HS} \cite[Theorem 1.3]{Simplicity}.  In other words, $\mathcal{M}_p$ is the one-dimensional manifold of extremal functions for \eqref{HS}.
It is known that $G_p$ is compact if $g\in L^q(\Omega),$ where $q=N/p$ when $p\in (1,N),$ $q>1$ when $p=N,$ and $q=1$ when $p>N$ \cite{cuesta,Nirjan}. For $p\in(1,N)$, this result was extended to $g\in L^{N/p,d}(\Omega)$ with $d<\infty$ \cite{Vis-05}, and further to the closure of $C_c^\infty(\Omega)$ in $L^{N/p,\infty}(\Omega)$ \cite{Anoop_Ejd}. Moreover, for $p\in (1,N)$, a necessary and sufficient condition for compactness is established in \cite{Anoo-_Das1}: $G_p$ is compact if and only if $g$ belongs to the closure of $C_c^\infty(\Omega)$ in $\H$.

\smallskip

 If $\Omega$ is bounded and $g=1$, then $G_p$ is compact, and \eqref{HS} reduces to the Poincar\'e inequality. The sharp gradient stability of the Poincar\'e inequality for bounded convex domains was proved in \cite[Theorem 3.1]
 {Amir2}. More precisely, there exists a constant $C>0$ such that
 \begin{equation*}\label{sharp-poincare}
     \frac{\|\nabla u\|_{L^p(\Omega)}^p-\mathbb{H}_p(\Omega)^p\|u\|_{L^p(\Omega)}^p}
{\|\nabla u\|_{L^p(\Omega)}^p}\geq Cd_p(u)^\alpha\quad \forall\,u\in\Dp\,,
 \end{equation*}
where $\alpha=\max\{2,p\}$, and this exponent is optimal. Here, $d_p(u)$ denotes the \textit{gradient-level distance} from $u$ to the manifold $\mathcal{M}_p$, given by
\begin{equation}\label{d-distance}
d_p(u):=\inf\left\{
\frac{\|\nabla u-\nabla v\|_{L^p(\Omega)}}
{\|\nabla u\|_{L^p(\Omega)}}:
v\in\mathcal{M}_p
\right\}.
\end{equation} In this article, we investigate quantitative stability of the weighted Hardy-type inequality \eqref{HS} for Hardy potentials $g\in\H$ for which the corresponding map $G_p$ is compact.
Define the
\textit{$p$-Hardy deficit functional} $\delta_p(u)$ by
$$
\delta_p(u):=
\frac{\|\nabla u\|_{L^p(\Omega)}^p-\mathbb{H}_p(\Omega)^p\|u\|_{L^p(\Omega,\,|g|)}^p}
{\|\nabla u\|_{L^p(\Omega)}^p}\qquad \text{for}\; u\in \Dp.
$$
Our first result shows that, for $p=2$, the gradient stability estimate $\delta_p(u)\geq Cd_p(u)$ holds with an optimal stability constant $C>0.$ Throughout this article, $\Phi$ denotes a non-negative extremal for \eqref{HS}.
\begin{thm}\label{Thm-p=2} 
Let $\Omega\subseteq\mathbb{R}^N$ be a domain (bounded if $N\leq 2$), and let $g\in \mathcal{H}_2(\Omega)$ be such that the map $G_2$ is compact. 
Then
    \begin{equation}\label{Eqn-sta-2}
       \delta_2(u)\ge \left(1-\frac{\mathbb{H}_2(\Omega)^2}{\Lambda(\Omega)^2}\right) d_2(u)^\alpha\qquad\forall\,u\in \Dt,
    \end{equation}  
    where $\alpha=2$ and $$ \Lambda(\Omega):=\inf\left\{\frac{\|\nabla u\|_{L^2(\Omega)}}{\| u\|_{L^2(\Omega,\,|g|)}}:u\in\Dt\setminus\{0\},\;\langle u,\Phi\rangle_{\mathcal{D}_0^{1,2}(\Omega)}=0 \right\}>\mathbb{H}_2(\Omega)\,.$$
Moreover, the constant $ \left(1-\frac{\mathbb{H}_2(\Omega)^2}{\Lambda(\Omega)^2}\right)$ is sharp and is attained by non-trivial functions. Furthermore, if $\Lambda(\Omega)$ is attained at $v\in\Dt\setminus\{0\},$ then the equality in \eqref{Eqn-sta-2} occurs at $u=c\,\Phi+v,$ for any $c\in \R$. Conversely, if equality in \eqref{Eqn-sta-2} holds at $u\in \Dt\setminus\{0\}$ then there exists $c\in \R$ such that $\Lambda(\Omega)$ is attained at $u+c\Phi.$
\end{thm}
\begin{rem}\label{rem---sharp=2}
The exponent $\alpha=2$ in \eqref{Eqn-sta-2} is sharp if $\Phi\in C_{\mathrm{loc}}^1(\Omega)$ (see Remark \ref{last-rem}). If $\Omega$ is bounded and $g\in L^q(\Omega)$ for some $q>2N$, then interior regularity gives $\Phi\in C_{\mathrm{loc}}^{1}(\Omega)$; see \cite[Corollary]{regulity}. Thus, $\alpha=2$ is sharp under this condition on $g.$ In particular, \eqref{Eqn-sta-2} gives sharp gradient stability estimate for the Poincar\'e inequality with the optimal stability constant.

\end{rem}

The proof of Theorem \ref{Thm-p=2} requires the standard Bianchi-Egnell-type argument in stability theory. But to show $\Lambda(\Omega)$ is attained and $\Lambda(\Omega)>\mathbb{H}_2(\Omega)$ we use the compactness of the map $G_2$ (see Proposition \ref{prop;compact}). To prove the sharpness and the attainment of the constant in \eqref{Eqn-sta-2}, we use the Hilbert space structure of the underlying spaces and exploit the fact that $\mathcal{M}_p$ is a one-dimensional manifold.

\smallskip

If $\Omega$ is a bounded domain, then, by the Poincar\'e inequality $\|\nabla u\|_{L^p(\Omega)}\geq C\|u\|_{L^p(\Omega)}$ for every $u\in \Dp,$ for some $C>0,$ we obtain $d_p(u)\geq CD_p(u)$, where $D_p(\cdot)$ denotes the distance function defined by
\begin{equation}\label{D-Distance}
D_p(u):=\inf\left\{\frac{\|u-v\|_{L^p(\Omega)}}{\|\nabla u\|_{L^p(\Omega)}}:v\in\mathcal{M}_p\right\}.
\end{equation}
Hence, as a corollary of Theorem \ref{Thm-p=2}, we obtain the following quantitative stability estimate: 
\begin{cor}
    Let $\Omega $ be a bounded domain and $g\in \mathcal{H}_2(\Omega)$ be such that the map $G_2$ is compact. 
Then there exists $C=(N,\Omega,g)>0$ such that
    \begin{equation*}
    \delta_2(u)\geq C\,D_2(u)^2\qquad\forall\,u\in \Dt.
    \end{equation*}  
\end{cor}

Next, we consider the nonlinear case $p\ne 2$ and assert that a quantitative stability estimate of the form $\delta_p(u)\ge C D_p(u)^\alpha$ holds for every $u\in \Dp,$ for some constants $C,\alpha>0,$ when the Hardy potential $g$ belongs to a certain class. 
\begin{thm}\label{Hardy-Thm-p} Let $\Omega\subset \R^N$ be a ball and $p\in (1,\infty)$. Let $g:\Omega\to \R$ be a radial function and
$g\in L^q(\Omega),$ where $q$ satisfies \begin{align} \label{weights}
   q  
\begin{cases} 
>\frac{N}{p}, & \text{if } p\in (1,N],\\
=1, & \text{if } p > N.
\end{cases} 
\end{align}
Then there exists $C=C(p,N,\Omega,g)>0$ such that 
\begin{equation}\label{Eqn-Hardy-p}
    \delta_p(u)\ge C D_p(u)^\alpha\,,\qquad\forall\,u\in \Dp\,,
\end{equation}
where $\alpha=\max\{2,p\}.$
\end{thm}

To prove Theorem \ref{Hardy-Thm-p}, we first establish a lower estimate for the deficit functional $\delta_p(u)$ in terms of the remainder 
$\|\Phi\nabla(u/\Phi)\|_{L^p(\Omega)}$
(see Lemma \ref{reminder-lemma}). Next, to estimate this remainder, we show that $\Phi$ is radially nonincreasing (see Proposition \ref{prop-radi-dec}). This allows us to apply a weighted Poincar\'e inequality for radially nonincreasing weights $\Phi^p$ (see Lemma \ref{Poincare}) to obtain $\|\Phi\nabla(u/\Phi)\|_{L^p(\Omega)}\gtrsim D_p(u).$ Finally, using these two estimates, we obtain the stability estimate \eqref{Eqn-Hardy-p}.
The assumptions on the Hardy potentials $g$ are used to establish both estimates.

\smallskip

Note that if $g\in L^\infty(\Omega)$, Theorem \ref{Hardy-Thm-p} yields the stability estimate \eqref{Eqn-Hardy-p}, which is weaker than $\delta_p(u)\gtrsim d_p(u)^\alpha$, since $d_p(\cdot)\gtrsim D_p(\cdot)$. We now prove this gradient stability estimate. 

\begin{thm}\label{thm-gard>2} Let $\Omega\subset\R^N$ be a ball and $p\in (1,\infty).$ Let $g:\Omega\to \R$ be a radial function and $g\in L^\infty(\Omega).$ Then there exists $C=(p,N,\Omega,g)>0$ such that
\begin{equation}\label{Eqn-Hardy-p-infty}
    \delta_p(u)\ge C d_p(u)^\alpha,\quad\forall\,u\in \Dp,\,
\end{equation}
where $\alpha=\max\{2,p\}.$ Moreover, the exponent $\alpha$ is sharp.
\end{thm}\begin{rem}
  The sharp gradient stability of the Poincar\'e inequality follows by setting $g=1$ in Theorem \ref{thm-gard>2}.
\end{rem}

\smallskip

We organize the article as follows: 
 Section \ref{main-results} is devoted to the proofs of the stability results (Theorems \ref{Thm-p=2}, \ref{Hardy-Thm-p}, and \ref{thm-gard>2}), as well as several auxiliary results used in their proofs. 

\smallskip

\textbf{Notation: } We use the following notation.

$(i)$ $f \gtrsim g$ means that there exists
a positive constant $C>0$ such that $f(x)\ge C g(x)$ for every $x.$ $(ii)$ $f \lesssim g$ means that there exists a
 positive constant $C>0$ such that $f(x)\le C g(x)$ for every $x.$ 

\section{Stability }\label{main-results} In this section, we study quantitative stability of the weighted Hardy-type inequality \eqref{HS}. We assume that $\Omega$ be a domain in $\R^N$ such that $\Omega\;\text{is bounded when}\;p\ge N,$ while $\Omega$ can be any domain if $p\in (1,N).$ Under this restriction on $p,N,\Omega$, it is known that the Beppo-Levi space $\Dp$ can be identified with a subspace of the space
of distributions \cite{Lions}. 
Recall, the best constant in \eqref{HS}
  \begin{equation}\label{best-const}
      \mathbb{H}_{p}(\Omega)=\inf\left\{R_p(u):u\in\Dp\setminus\{0\}\right\},\quad R_p(u)=\frac{\|\nabla u\|_{L^p(\Omega)}}{\|u\|_{L^p(\Omega,\,|g|)}}\,.
  \end{equation} 
If the best constant $\mathbb{H}_{p}(\Omega)$ is attained at some $\Phi\in \D_0^{1,p}(\Omega)\setminus\{0\}$, then $\Phi$ is an extremal for \eqref{HS}. Moreover, an extremal (if exists) satisfies the Euler-Lagrange equation 
\begin{equation}\label{Euler-Lag-p}
    \int_\Omega |\nabla \Phi|^{p-2}\nabla \cdot \nabla v\dd x=\mathbb{H}_{p}(\Omega)^p\,\int_\Omega |g||\Phi|^{p-2}\,\Phi\,v\dd x,\quad\forall\, v\in \Dp.
\end{equation}
In other words, $\Phi$ solves the following nonlinear partial differential equation weakly:
\begin{equation*}\label{P} \tag{$\mathcal{P}$}   -\Delta_p u=\mathbb{H}_{p}(\Omega)^p\,|g||u|^{p-2}u\quad\text{in}\;\Dp,
\end{equation*}
where $\Delta_p u=\text{div}(|\nabla u|^{p-2}\nabla u).$

Our goal is to investigate quantitative stability of the weighted Hardy-type inequality \eqref{HS} when $\mathbb{H}_p(\Omega) $ is attained. Of course, there are Hardy potentials $g\in \H$ for which $\mathbb{H}_p(\Omega) $ is not attained, i.e., $\mathcal{M}_p=\{0\}$.  
A sufficient condition for the existence of non-trivial extremals for \eqref{HS} is the compactness of the map $G_p$. It is well known that if  $\mathbb{H}_p(\Omega)$ is attained then $\mathbb{H}_p(\Omega)^p$ is the first (i.e., the least) weighted eigenvalue of the problem \eqref{P}, and it is also principal (i.e., the corresponding eigenfunctions do not change their sign in $\Omega$) 
\cite[Theorem 1.3]{Simplicity}. Moreover, the weighted eigenvalue $\mathbb{H}_p(\Omega)^p$ is simple (i.e.,  the corresponding eigenfunctions are constant multiples of each other). Thus, for $g\in\H$ if $G_p$ is compact, then the extremal functions in \eqref{HS} are given by $$ \mathcal{M}_p=\{c\,\Phi: c\in \mathbb{R}\},$$ where $\Phi\in \Dp\setminus\{0\}$ is a non-negative solution of \eqref{P}. 

Moreover, if $g\in L^q(\Omega),$ where $q$ is given by \eqref{weights}, and $\Omega$ is bounded, then $\Phi$ is locally H\"older continuous \cite[Theorem 8]{Serrin} and $\Phi>0$ in $\Omega$ \cite[Proposition 3.2]{cuesta}. Furthermore, $\Phi\in C_{\mathrm{loc}}^{1}(\Omega)$ when $g\in L^r(\Omega)$ for some $r>Np'$ \cite[Corollary]{regulity}.

\subsection{The linear case:} $p=2$. First, we study the gradient stability of the weighted Hardy-type inequality \eqref{HS} in the linear case $p=2.$ Recall that $$G_2(u)=\int_\Omega|g||u|^2\dd x\,\quad\text{and}\quad \Lambda(\Omega)=\inf\left\{R_2(u):u\in\Dt\setminus\{0\},\;\langle u,\Phi\rangle_{\Dt}=0 \right\}.$$ 
The following proposition shows that if $G_2$ is compact, then $\Lambda(\Omega)$ is attained and strictly larger than the Hardy constant $\mathcal{H}_2(\Omega)$. 
\begin{prop}\label{prop;compact} Let $\Omega\subseteq\mathbb{R}^N$ be a domain (bounded if $N\leq 2$), and let $g\in \mathcal{H}_2(\Omega)$ be such that the map $G_2$ is compact.  Then $\Lambda(\Omega)$ is attained in $\Dt\setminus\{0\}$ and $$\Lambda(\Omega)>\mathbb{H}_2(\Omega).$$
\end{prop}
\begin{proof}
By the definition of $\Lambda(\Omega),$  
there exists $(v_k)$ in $\D_0^{1,2}(\Omega)$ with $\langle v_k,\Phi\rangle_{\Dt}=0$ such that $$R_2(v_k)\to \Lambda(\Omega).$$ Since the Rayleigh quotient $R_2(\cdot)$ is homogeneous, by normalizing the minimizing sequence (i.e., $\|v_k\|_{L^2(\Omega,\,|g|)}=1$) and calling it again $(v_k)$, we get $\|\nabla v_k\|_{L^2(\Omega)}\to \Lambda(\Omega).$ Thus, $(v_k)$ is a bounded sequence in $\D_0^{1,2}(\Omega),$ and hence it has a weakly convergent subsequence, call it again $v_k\rightharpoonup v_0$ in $\Dt$. Therefore, $$\int_\Omega\nabla \Phi\cdot\nabla v_0\dd x=\lim_{k\to\infty}\int_\Omega\nabla \Phi\cdot\nabla v_k\dd x=0,$$ i.e., $\langle v_0,\Phi\rangle_{\Dt}=0.$
Moreover, weak lower semi-continuity gives $$\int_\Omega|\nabla v_0|^2\,\dd x\le \lim_{k\to \infty}\inf\int_\Omega|\nabla v_k|^2\,\dd x=\Lambda(\Omega)^2.$$ Since $G_2$ is compact, we infer that $$\int_\Omega|g||v_0|^2\dd x=\lim_{k\to \infty}\int_\Omega|g||v_k|^2\dd x=1.$$ Consequently, $R_2(v_0)\le \Lambda(\Omega),$ and hence $R_2(v_0)= \Lambda(\Omega),$ i.e., $\Lambda(\Omega)$ is attained at $v_0.$

Now we prove $\Lambda(\Omega)> \mathbb{H}_2(\Omega).$ Clearly $\Lambda(\Omega)\ge \mathbb{H}_2(\Omega).$  If we assume, $\Lambda(\Omega)=\mathbb{H}_{2}(\Omega),$ then $v_0$ is an extremal for \eqref{HS}. Therefore, $v_0=a_0\Phi$ for some $a_0\ne 0,$ contradicts $\langle v_0,\Phi\rangle_{\Dt}=0.$ Hence the proof.
\end{proof}
Now we use Proposition \ref{prop;compact} to prove the gradient stability of \eqref{HS} in the linear case $p=2.$
\begin{proof}[Proof of Theorem \ref{Thm-p=2}]
It it enough to prove the stability inequality \eqref{Eqn-sta-2} for all $u\in\Dt$ satisfying $$\|\nabla u\|_{L^2(\Omega)}=1.$$ Choose such $u$ arbitrary.
Since $\mathcal{M}_2=\{c\Phi:c\in\R\}$ is a one-dimensional subspace of the Hilbert space $\Dt$, it is closed. Thus, by the Hilbert projection theorem,  $d_2(u)$ is uniquely attained \cite[Theorem 7.1.2, Remark 7.1.8]{kesh-f}, i.e., there exists a unique $a_0=a_0(u)\in\mathbb{R}$ such that 
$$d_2(u)=\| \nabla u- a_0\nabla\Phi\|_{L^2(\Omega)} \,.$$
Moreover, we can decompose $u$ as 
\begin{equation}\label{decompose}
    u=a_0\Phi+d_2(u)
\,v,
\end{equation}
 where  $\langle v,\Phi\rangle_{\Dt}=0$ with $\| \nabla v\|_{L^2(\Omega)}=1$. Assume that $\|\nabla\Phi\|_{L^2(\Omega)}=1.$
Now, we have
\begin{align}
\int_\Omega|\nabla u|^2{\mathrm{d}}x&=a_0^2\int_\Omega |\nabla \Phi|^2\,{\mathrm{d}}x+2a_0\,d_2(u)\int_\Omega\nabla \Phi\cdot\nabla v\,{\mathrm{d}}x+d_2(u)^2\int_\Omega|\nabla v|^2{\mathrm{d}}x\nonumber    \\& =a_0^2+d_2(u)^2.\label{p=2-1-1}
\end{align}
On the other hand, we have
\begin{align}
\int_\Omega|g||u|^2{\mathrm{d}}x&=a_0^2\int_\Omega|g||\Phi|^2{\mathrm{d}}x+2a_0\,d_2(u)\int_\Omega|g|\Phi\, v\,{\mathrm{d}}x+d_2(u)^2\int_\Omega|g||v|^2{\mathrm{d}}x\nonumber\\&= \frac{a_0^2}{\mathbb{H}_2(\Omega)^2}+2a_0d_2(u)\int_\Omega|g|\Phi \,v\,{\mathrm{d}}x+d_2(u)^2\int_\Omega|g||v|^2{\mathrm{d}}x.  \label{p=2--2}
\end{align}
Since $\langle v,\Phi\rangle_{\Dt}=0$ and $-\Delta \Phi=\mathbb{H}_2(\Omega)^2|g|\Phi$ (see \eqref{P}), it follows that $$\int_\Omega|g|\Phi\,v\,{\mathrm{d}}x =\frac{1}{\mathbb{H}_2(\Omega)^2}\int_\Omega \nabla \Phi\cdot\nabla v\,{\mathrm{d}}x=0.$$ Thus, from \eqref{p=2--2}, we get $$\mathbb{H}_2(\Omega)^2\int_\Omega|g||u|^2{\mathrm{d}}x=a_0^2+d_2(u)^2\mathbb{H}_2(\Omega)^2\int_\Omega|g||v|^2{\mathrm{d}}x,$$ combining this with \eqref{p=2-1-1}, we obtain $$
    \delta_2(u)=\int_\Omega|\nabla u|^2{\mathrm{d}}x- \mathbb{H}_2(\Omega)^2\int_\Omega|g||u|^2\,{\mathrm{d}}x= d_2(u)^2\left(1-\mathbb{H}_2(\Omega)^2\int_\Omega|g||v|^2{\mathrm{d}}x\right).
$$
Due to $\langle v,\Phi\rangle_{\Dt}=0$ and $\| \nabla v\|_{L^2(\Omega)}=1,$ 
$$ \Lambda(\Omega)^2\le \frac{\int_\Om |\nabla v|^2{\mathrm{d}}x}{\int_\Omega|g||v|^2{\mathrm{d}}x}= \frac{1}{\int_\Omega|g||v|^2{\mathrm{d}}x},$$ implying that $\int_\Omega|g||v|^2{\mathrm{d}}x\le \frac{1}{\Lambda(\Omega)^2}.$ Therefore, $$\delta_2(u)\ge d_2(u)^2\left(1-\frac{\mathbb{H}_2(\Omega)^2}{\Lambda(\Omega)^2}\right),$$ proving the stability inequality \eqref{Eqn-sta-2}. Now we prove that equality in \eqref{Eqn-sta-2} is attained in $\Dt\setminus\{0\}$.

We know that $\Lambda(\Omega)$ is attained, see Proposition \ref{prop;compact}.  Assume that $\Lambda(\Omega)$ is attained at $v\in\Dt\setminus{0}$.
Take $u\in \Dt$ of the form $u=c\Phi+v$ for some arbitrary $c\in \R.$ Since $\langle v,\Phi\rangle_{\Dt}=0,$
 a similar computation as before gives
\begin{align*}
\delta_2(u)&=\frac{1}{\|u\|^2_{L^2(\Omega,\,|g|)}}\Big[c^2\int_\Omega|\nabla \Phi|^2\,{\mathrm{d}}x+\int_\Omega|\nabla v|^2\,{\mathrm{d}}x\\&\qquad-c^2\,\mathbb{H}_2(\Omega)^2\int_\Omega|g||\Phi|^2{\mathrm{d}}x-\mathbb{H}_2(\Omega)^2\int_\Omega|g||v|^2{\mathrm{d}}x\Big]\\&=\frac{1}{\|u\|^2_{L^2(\Omega,\,|g|)}}\Big[\|\nabla v\|^2_{L^2(\Omega)}-\mathbb{H}_2(\Omega)^2\|v\|^2_{L^2(\Omega,\,|g|)}\Big]=\left(1-\frac{\mathbb{H}_2(\Omega)^2}{\Lambda(\Omega)^2}\right)\frac{\|\nabla v\|^2_{L^2(\Omega)}}{\|u\|^2_{L^2(\Omega,\,|g|)}}\,,
\end{align*}
where we used the attainability of $\Lambda(\Omega)$  at $v$ in the last equality.
On the other hand,
 \begin{align*}  d_2(u)^2=&\frac{1}{\|u\|^2_{L^2(\Omega,\,|g|)}}\inf_{c\in\mathbb{R}}\int_\Omega|(a_0-c)\nabla \Phi+\nabla v|^2\,{\mathrm{d}}x\\=&\frac{1}{\|u\|^2_{L^2(\Omega,\,|g|)}}\Big[\inf_{c\in\mathbb{R}}(a_0-c)^2\int_\Omega|\nabla \Phi|^2\,{\mathrm{d}}x+\int_\Omega|\nabla v|^2\,{\mathrm{d}}x\Big]=\frac{\|\nabla v\|^2_{L^2(\Omega)}}{\|u\|^2_{L^2(\Omega,\,|g|)}}\,,
 \end{align*}
 which shows that the equality in \eqref{Eqn-sta-2} is attained at $u\ne 0$ and hence the constant $ \left(1-\frac{\mathbb{H}_2(\Omega)^2}{\Lambda(\Omega)^2}\right)$ is sharp. 
 
Conversely, assume that the equality occurs in \eqref{Eqn-sta-2} for some $u\in \Dt\setminus\{0\}.$ Notice that $u$ can be decomposed as $u=a_0\Phi+v,$ for some $a_0=a_0(u)\in \R$ and $\langle v,\Phi\rangle_{\Dt}=0$ (see \eqref{decompose}). Then we have as before
 $$
\delta_2(u)=\frac{1}{\|u\|^2_{L^2(\Omega,\,|g|)}}\Big[\|\nabla v\|^2_{L^2(\Omega)}-\mathbb{H}_2(\Omega)^2\|v\|^2_{L^2(\Omega,\,|g|)}\Big]\quad\text{and}\quad d_2(u)^2=\frac{\|\nabla v\|^2_{L^2(\Omega)}}{\|u\|^2_{L^2(\Omega,\,|g|)}}.  $$ Since $\delta_2(u)= \left(1-\frac{\mathbb{H}_2(\Omega)^2}{\Lambda(\Omega)^2}\right)d_2(u)^2,$ it follows that
 $$\|\nabla v\|_{L^2(\Omega)}=\Lambda(\Omega)\|v\|_{L^2(\Omega,\,|g|)}  .$$ That is $\Lambda(\Omega)$ is attained at $v.$ This completes the proof.
\end{proof}

\subsection{Non-linear case:} $p\ne 2.$ In this subsection, we investigate quantitative stability of the weighted Hardy-type inequality \eqref{HS} for $p\neq 2$. 
We recall a weighted Poincar\'e inequality which will be used to prove stability. The proof follows from \cite[Corollary 3]{Dyda_poin}. 
 
 We call a function $\phi : B_R(x_0) \to [0, \infty)$ a radially nonincreasing weight if $\phi$ is a radial function, i.e.
$\phi(x)=\Psi(|x-x_0|)$ and its profile $\Psi$ is nonincreasing and right-continuous with left-limits. 
\begin{lem}\label{Poincare}
    Let $p\ge 1$ and $\phi: B_R(x_0) \mapsto[0,\infty)$ be a radially nonincreasing function. There exists a positive
constant $C = C(N,p,R,\phi)$ such that
\begin{equation*}
   \int_{B_R(x_0)}|\nabla u|^p\phi\,{\mathrm{d}}x\ge C \int_{B_R(x_0)}|u-(u)_\phi|^p\phi\,{\mathrm{d}}x\,,
\end{equation*}
for every $u\in W^{1,p}(B_R(x_0)),$ where $(u)_\phi=\left(\int_{B_R(x_0)} \phi\, {\mathrm{d}}x\right)^{-1}\int_{B_R(x_0)} u \,\phi \,{\mathrm{d}}x.$
\end{lem}
Next, we recall a lemma describing some properties of the extremals of \eqref{HS}.
\begin{lem}\cite[Proposition 2.1, Proposition 3.2]{cuesta}\cite[Corollary]{regulity}\label{cuesta-lem} Let $\Omega\subset \mathbb{R}^N$ be a bounded domain and $g\in L^q(\Omega)$, where $q$ is given by \eqref{weights}. If $u$ is a non-negative solution of \eqref{P}, then $u\in L^\infty(\Omega)$, is locally H\"older continuous, and is positive in $\Omega$. Moreover, if  $g\in L^r(\Omega)$ for some $r>Np'$, then $u\in C_{\mathrm{loc}}^{1,\gamma}(\Omega)$ for some $\gamma\in (0,1)$.
\end{lem}
To use Lemma \ref{Poincare} to prove quantitative stability of \eqref{HS}, we need to show that its non-negative extremals are radially nonincreasing. We prove this under the assumption that the Hardy potentials $g\in\H$ is radial. The proof is also of independent interest. 
\begin{prop}\label{prop-radi-dec}
 Let $\Omega=B_R(x_0)$ and $p\in (1,\infty).$  Let $g:\Omega\to \R$ be a radial function such that $g\in L^q(\Omega),$ where $q$ is given by \eqref{weights}. Then any non-negative extremal of \eqref{HS} is radial and radially nonincreasing.
\end{prop}
\begin{proof} By translation, we may assume, without loss of generality, that 
 $\Omega=B_R(0).$ Let $u\ne 0$ be a non-negative extremal of the weighted Hardy-type inequality \eqref{HS}. By Lemma \ref{cuesta-lem}, we have $u\in L^\infty(\Omega)$ and $u$ is locally H\"older continuous in $\Omega$. First, we show that $u$ is radial.
    Let $Q\in S(N)$ be an orthogonal rotation and define $$v(x):=u(Qx)\quad \text{for}\;x\in \Omega.$$ Since $Q$ is orthogonal $\nabla v(x)=Q^T\nabla u(Qx)$ and consequently, $|\nabla v(x)|=|\nabla u(Qx)|
    .$ Thus, a standard computation shows that the $p$-Laplacian is invariant under orthogonal changes of variables, and hence $$\Delta_p v(x)=(\Delta_p u)(Qx)\,,\quad x\in \Omega.$$ Since $u$ solves \eqref{P} and $g$ is radial,  
    \begin{align*}
   -\Delta_p v(x)=\mathbb{H}_p(\Omega)^p\,|g(|Qx|)||u(Qx)|^{p-2}u(Qx)=\mathbb{H}_p(\Omega)^p\,|g(|x|)||v(x)|^{p-2}v(x).          
    \end{align*}
  Moreover, $v=0$ on $\partial\Omega,$ and hence $v$ is also an extremal for \eqref{HS}. Thus, there exists a constant $C>0$ such that $v=Cu.$  Using the rotation preserves the normization and $g$ is radial $$\int_\Omega |g|\,|v|^p\dd x=\int_\Omega |g|\,|u|^p{\mathrm{d}}x>0,$$ it follows that $C^p=1,$ and hence $u(x)=v(x)=u(Qx).$ Since $Q$ is an arbitrary orthogonal rotation, $u$ is a radial function.

   Now, we prove that $u$ is nonincreasing.
Since $u$ satisfies \eqref{P}, it can be verified that $u(x)=U(|x|)$ satisfies  
\begin{equation}\label{el-radi}
             -\left(r^{N-1}|U'|^{p-2}U'\right)' = \mathbb{H}_p(\Omega)^p\, r^{N-1}|g|U^{p-1} \qquad \text{in }\mathcal{D}(0,R)\,,
            \end{equation}
where $\mathcal{D}(0,R)$ is the space of distributions. Define $$ F(r):=r^{N-1}|U'(r)|^{p-2}U'(r)\quad \text{and}\quad Q(r):=\mathbb{H}_p(\Omega)^p  r^{N-1}|g(r)|
U^{p-1}(r) . $$ 
Then \eqref{el-radi} can be rewrite as $F'=-Q$ in $\mathcal{D}(0,R).$  Observe that 
$$\int_{\mathbb{S}^{N-1}} \dd \theta \int_0^R|U'(r)|^{p-1}r^{N-1} \dd r= \int_\Omega |\nabla u|^{p-1}\dd x\le\left(\int_\Omega\dd x\right)^\frac{1}{p}\left(\int_\Omega|\nabla u|^p\dd x\right)^\frac{p-1}{p}<\infty\,, $$ which implies that
$F\in L^1((0,R)).$ Moreover, using the fact that $u\in L^\infty(\Omega)$ and $g\in L^1(\Omega),$ we get $Q\in L^1((0,R)).$ Consequently, the distribution derivatives of $F$ also in $L^1((0,R))$, and hence $F\in W^{1,1}((0,R)).$ Which implies that $F$ is absolutely continuous (see \cite[Theorem 2.1.3]{Kesh_PDE}). Using this, 
 next, we prove $F(0)=0.$
 Towards this, we consider test functions $v_\epsilon(x)=\psi_\epsilon(|x|)\in \Dp,$ where
\[ \psi_\epsilon(r)=  
\begin{cases} 
1 & \text{if } 0<r\le \epsilon,\\
\frac{2\epsilon-r}{\epsilon} & \text{if }\epsilon< r\le  2\epsilon,\\
0 & \text{if }2\epsilon< r< R.
\end{cases}\qquad 0<\epsilon< R/2\]
Use this as a test function in the Euler Lagrange equation \eqref{Euler-Lag-p}, we get $$\int_0^RF(r)\psi_\epsilon'(r)\dd r= \int_0^RQ(r)\psi_\epsilon(r) \dd r= \int_0^{2\epsilon}Q(r)\psi_\epsilon(r)  \dd r\le \int_0^{2\epsilon}Q(r)  \dd r.$$ 
Therefore, 
\begin{equation}\label{eq1:avg}
\left|\frac{1}{\epsilon}\int_\epsilon^{2\epsilon}F(r)\right| =\left|\int_0^RF(r)\psi_\epsilon'(r)\dd r\right|\to 0\quad{\text{as}}\quad \epsilon\to 0 .    
\end{equation}
On the other hand, since $F'=-Q$ in $\mathcal{D}(0,R)$ and $F$ is absolutely continuous, 
 $$|F(r)-F(\epsilon)|\le \int_\epsilon^{2\epsilon} Q(s)\dd s,\quad \forall\,r\in [\epsilon,2\epsilon].$$ Averaging this estimate gives
 \begin{equation}\label{eq2:avg}
 \left|F(\epsilon)-\frac{1}{\epsilon}\int_\epsilon^{2\epsilon}F(r)\dd r\right|\le\int_\epsilon^{2\epsilon} Q(s)\dd s\to 0 \quad\text{as}\quad \to 0 .  
 \end{equation}
Combining \eqref{eq1:avg} and \eqref{eq2:avg} we conclude that $F(0)=0.$ Thus, the fundamental theorem of calculus for absolutely continuous functions gives
  $$F(r) = \int_0^r F'(s)\,\dd s=-\int_0^r Q(s)\,\dd s= 
-\int_0^r \mathbb{H}_p(\Omega)^p  s^{N-1}|g(s)|U^{p-1}(s)\,\dd s\le 0.$$ Which implies that  $U'(r)\le0 $ for a.e. on $(0,R).$ Since $U\in W^{1,p}(a,b)$ for all $[a,b]\subset (0,R),$ so $U$ is absolutely continuous on $[a,b],$ and $U(b)-U(a)=\int_a^bU'(r)\dd r\le 0.$ Therefore, $U$ is a nonincreasing function on $(0,R)$, completing the proof. 
        \end{proof}
\begin{rem}
 If $g$ is radially decreasing, Proposition \ref{prop-radi-dec} was proved in \cite[Theorem 1.2]{Anoop_Ejd} for $\Omega=B_R(0)$ using Schwarz symmetrization. Thus, Proposition \ref{prop-radi-dec} is more general.
\end{rem}

 We recall a lemma, which will be used to prove  quantitative stability of the weighted Hardy-type inequality \eqref{HS} (Theorems \ref{Hardy-Thm-p}, \ref{thm-gard>2}). 
This lemma concerns the equivalence of the {\it energy functional} and the {\it simplified energy functional}. We will replace the role of the following energy functional associated to \eqref{P}
$$Q_{p}(u)=\int_\Omega[|\nabla u|^p- \mathbb{H}_p(\Omega)^p|g||u|^p] \ {\mathrm{d}}x\quad\,\forall\, u\in W^{1,p}(\Omega) \cap C_c(\Omega)$$
by the simplified energy functional given by
$$Q_{p}^{sim}(u)=\int_\Omega w^2|\nabla u|^2\left(u|\nabla w|+w|\nabla u|\right)^{p-2}{\mathrm{d}}x \quad\,\forall\, u\in W^{1,p}(\Omega) \cap C_c(\Omega) \,,$$
where $w  \in W^{1,p}_\mathrm{loc}(\Omega)\cap C(\Omega)$ is a positive solution of the following quasi-linear elliptic equation:
\begin{equation*}
    \mathcal{Q}'_{p}[\phi]:=-\Delta_p \phi- \mathbb{H}_p(\Omega)^p|g| |\phi|^{p-2}\phi=0 \ \,\text{in} \, \ \Omega \,.
\end{equation*}
\begin{lem}\cite[Lemma 2.2]{Yehuda2008} and \cite[Lemma 3.4]{Yehuda2015}\label{prop-1}
    Let $p\in (1,\infty)$ and $w\in W^{1,p}_\mathrm{loc}(\Omega)\cap C(\Omega)$ be a positive solution of $\mathcal{Q}_{p}^\prime[\phi]=0$. Then $\exists\,c,\,C>0$ depending only on $p$ such that $$c\,Q_{p}( u)\le  Q_{p}^{sim}(\frac{u}{w})\le C\,Q_{p}( u)\quad  $$
    for all non-negative $u\in W^{1,p}(\Omega) \cap C_c(\Omega)$.
\end{lem}
The following lemma proves an estimate of the deficit functional, which will be used to prove Theorem \ref{Hardy-Thm-p} and Theorem \ref{thm-gard>2}.  
\begin{lem}\label{reminder-lemma}
    Let $\Omega\subset \R^N$ be a bounded domain and $p\in (1,\infty)$. Let 
$g\in L^q(\Omega),$ where $q$ satisfies \eqref{weights}. Then, there exists $C=C(p)>0$ such that 
\begin{align} \label{rem...}
\int_\Omega|\nabla u|^p\dd x-\mathbb{H}_p(\Omega)^p\int_\Omega|g||u|^p\dd x\ge C 
\begin{cases}
[I(u,\Phi)]^p
& \text{if } p\ge 2, \\[1ex]
\min\left\{[I(u,\Phi)]^p,\,\left[ I(u,\Phi)\right]^2\right\}\
& \text{if}\; p\in (1,2)
\end{cases}
\end{align} 
holds for every $u\in \Dp$ satisfying $\|\nabla u\|_{L^p(\Omega)}=1,$ where $$I(u,\Phi)=\left(\int_\Omega \left|\nabla \left(\frac{u}{\Phi}\right) \right|^p\Phi^p\,{\mathrm{d}}x\right)^{1/p}.$$
\end{lem}
\begin{proof} Since $ C_c^1(\Omega)$ is dense in $\Dp,$ it is sufficient to prove \eqref{rem...}
for all $u\in  C_c^1(\Omega)$ satisfying $\|\nabla u\|_{L^p(\Omega)}=1.$ Choose such $u$ arbitrary. Notice that
    $$Q_{p}(u)=\int_\Omega|\nabla u|^p\dd x-\mathbb{H}_p(\Omega)^p\int_\Omega|g||u|^p\dd x=Q_{p}(|u|)\,.
$$ 
It is known that $\Phi$ is a positive continuous solution of $\mathcal{Q}_{p}^\prime[\phi]=0$ (see Lemma \ref{cuesta-lem}). Thus, by Lemma \ref{prop-1} we obtain 
\begin{align}
  Q_p(u)  & \ge C Q_p^{sim}(\frac{|u|}{\Phi}) =C\int_\Omega \Phi^2\left|\nabla \left(\frac{|u|}{\Phi}\right) \right|^2\left(\frac{|u|}{\Phi}|\nabla \Phi|+\Phi\left|\nabla \left(\frac{|u|}{\Phi}\right) \right|\right)^{p-2}{\mathrm{d}}x \label{p2...} \,,
\end{align}
for some $C=C(p)>0.$ 

Case 1 : $p\ge 2.$ Using $p- 2\ge 0$ in \eqref{p2...} we obtain
\begin{equation*}
Q_p(u)   \ge C\int_\Omega \Phi^p\left|\nabla \left(\frac{|u|}{\Phi}\right) \right|^p{\mathrm{d}}x =C \,[I(u,\Phi)]^p.
\end{equation*}

Case 2 : $p\in (1,2).$ 
Set $\tilde{u}=\frac{|u|}{\Phi}.$ 
One can verify that for $p\in(1,2)$ $$(|a|+|b|)^{p-2}\ge \min\{|a|^{p-2},|b|^{p-2}\},\quad a,b\in \mathbb{R}.$$
Use this with $a=\tilde{u}\nabla \Phi$ and $b=\Phi\nabla \tilde{u}$ on the right hand side of \eqref{p2...} to get
\begin{align}
    Q_p(u) &\ge C \int_\Omega \Phi^2|\nabla \tilde{u}|^2\min\left\{ |\tilde{u}|^{p-2}|\nabla \Phi|^{p-2},\;\Phi^{p-2}|\nabla \tilde{u}|^{p-2}\right\}{\mathrm{d}}x\nonumber\\&= C\int_{A^c\cap \Omega} \Phi^p|\nabla \tilde{u}|^p {\mathrm{d}}x+ \int_{A}\Phi^2|\nabla \tilde{u}|^2|\nabla \Phi|^{p-2}|\tilde{u}|^{p-2}{\mathrm{d}}x\,,\label{eq1:p=2}
\end{align}
where $$A:=\{x\in \Omega: \Phi^{p-2}|\nabla \tilde{u}|^{p-2}>|\tilde{u}|^{p-2}|\nabla \Phi|^{p-2}\}.$$ Now, using H\"older inequality for conjugate pair $\left(\frac{2}{p},\frac{2}{2-p}\right)$ with $f=(\Phi|\nabla \tilde{u}|)^p(|\tilde{u}||\nabla \Phi|)^\frac{p(p-2)}{2}$ and $g=(|\nabla \Phi||\tilde{u}|)^\frac{p(2-p)}{2},$ we obtain 
\begin{equation}\label{Eq2:p=2}
    \int_A\Phi^p|\nabla \tilde{u}|^p{\mathrm{d}}x \le \left(\int_{A}\Phi^2|\nabla \tilde{u}|^2|\nabla \Phi|^{p-2}|\tilde{u}|^{p-2}{\mathrm{d}}x\right)^\frac{p}{2}\left(\int_{A}|\nabla \Phi|^p|\tilde{u}|^{p}{\mathrm{d}}x\right)^\frac{2-p}{2}.
\end{equation}
Combining \eqref{eq1:p=2} and \eqref{Eq2:p=2},  we deduce that
\begin{align} 
Q_p(u)\ge &C \int_{A^c\cap \Omega}\Phi^p|\nabla \tilde{u}|^p{\mathrm{d}}x+ \left(\int_A\Phi^p|\nabla \tilde{u}|^p{\mathrm{d}}x\right)^\frac{2}{p}\left(\int_{A}|\nabla \Phi|^p|\tilde{u}|^{p}{\mathrm{d}}x\right)^\frac{p-2}{p}\nonumber \\&\ge C \int_{A^c\cap \Omega}\Phi^p|\nabla \tilde{u}|^p{\mathrm{d}}x+ \left(\int_A\Phi^p|\nabla \tilde{u}|^p{\mathrm{d}}x\right)^\frac{2}{p}\left(\int_{\Omega}|\nabla \Phi|^p|\tilde{u}|^{p}{\mathrm{d}}x\right)^\frac{p-2}{p}\label{Eq6<p}.
\end{align}
where we used $p-2<0$ is the last inequality. 
Next, we prove the following estimate: \begin{align} \label{estT}
    \| \tilde{u}\nabla\Phi\|_{L^p(\Omega)}\le  \left(\frac{p}{p-1} \right) \|\nabla u\|_{L^p(\Omega)}\,.
\end{align} To see this, take $\varphi_\epsilon=\frac{|u|^p}{(\Phi+\epsilon)^{p-1}},\,\epsilon>0$ as a test function in weak formulation \ref{Euler-Lag-p} to get 
\begin{equation}\label{EVP1<2}
    \int_\Omega|\nabla \Phi|^{p-2}\nabla \Phi\cdot\nabla \varphi_\epsilon\,{\mathrm{d}}x=\mathbb{H}_p(\Omega)^p\int_\Omega |g|\Phi^{p-1}\,\varphi_\epsilon\,{\mathrm{d}}x\,.
\end{equation}
Since $$\nabla \varphi_\epsilon=p\frac{|u|^{p-2}u\nabla u}{(\Phi+\epsilon)^{p-1}}-(p-1)\frac{|u|^{p}\nabla \Phi}{(\Phi+\epsilon)^{p}}$$ using this in \eqref{EVP1<2} and rearranging we obtain
\begin{multline*}
     (p-1) \int_\Omega|\nabla \Phi|^{p}\left|\frac{u}{\Phi+\epsilon}\right|^p\,{\mathrm{d}}x=p\int_\Omega\frac{|\nabla\Phi|^{p-2}|u|^{p-2}u\nabla \Phi\cdot \nabla u }{(\Phi+\epsilon)^{p-1}}\,{\mathrm{d}}x\\-\mathbb{H}_p(\Omega)^p\int_\Omega|g|\left(\frac{\Phi}{\Phi+\epsilon} \right)^{p-1}|u|^p\,{\mathrm{d}}x.
\end{multline*}
Moreover, use H\"older inequality to obtain 
\begin{equation*}
    \int_\Omega\frac{|\nabla\Phi|^{p-2}\nabla \Phi\cdot |u|^{p-2}u\,\nabla u}{(\Phi+\epsilon)^{p-1}}\,{\mathrm{d}}x\le \left(\int_\Omega|\nabla \Phi|^{p}\left|\frac{u}{\Phi+\epsilon}\right|^p\,{\mathrm{d}}x\right)^\frac{p-1}{p}\left(\int_\Omega|\nabla u|^p\,{\mathrm{d}}x\right)^\frac{1}{p}\,.
\end{equation*}
Consequently, 
$$ (p-1) \int_\Omega|\nabla \Phi|^{p}\left|\frac{u}{\Phi+\epsilon}\right|^p\,{\mathrm{d}}x\le p \left(\int_\Omega|\nabla \Phi|^{p}\left|\frac{u}{\Phi+\epsilon}\right|^p\,{\mathrm{d}}x\right)^\frac{p-1}{p}\left(\int_\Omega|\nabla u|^p\,{\mathrm{d}}x\right)^\frac{1}{p},$$ and from this we infer $$\int_\Omega|\nabla \Phi|^{p}\left|\frac{u}{\Phi+\epsilon}\right|^p\,{\mathrm{d}}x\le\left(\frac{p}{p-1}\right)^p \int_\Omega|\nabla u|^p\,{\mathrm{d}}x.$$ Hence, the Monotone convergence theorem gives the desired estimate \eqref{estT}. Thus using this estimate in \eqref{Eq6<p}, together with the facts that $p-2<0$ and $\|\nabla u\|_{L^p(\Omega)}=1,$ we obtain 
\begin{equation}\label{equation:p<2:stability..}
 Q_p(u) \ge C_2\left[ \int_{A^c\cap \Omega}\Phi^p|\nabla \tilde{u}|^p{\mathrm{d}}x+ \left(\int_A\Phi^p|\nabla \tilde{u}|^p{\mathrm{d}}x\right)^\frac{2}{p}\right],
\end{equation}
for some $C_2=C_2(p)>0.$
Using the following elementary inequality $$a+b^\frac{2}{p}\ge C_p\min\{(a+b),(a+b)^\frac{2}{p}\}\quad \  \forall a,b\ge0$$
where $C_p=C_p(p)>0$, the right hand side of \eqref{equation:p<2:stability..} estimate as
$$\int_{A^c\cap \Omega}\Phi^p|\nabla \tilde{u}|^p{\mathrm{d}}x+ \left(\int_A\Phi^p|\nabla \tilde{u}|^p{\mathrm{d}}x\right)^\frac{2}{p}\ge C_p\min\left\{\int_{\Omega}\Phi^p|\nabla \tilde{u}|^p{\mathrm{d}}x,\left( \int_{\Omega}\Phi^p|\nabla \tilde{u}|^p{\mathrm{d}}x\right)^\frac{2}{p}\right\}.$$ Consequently,
\begin{equation*}
   Q_p(u)\ge C_2  C_p\min\left\{\int_{\Omega}\Phi^p|\nabla \tilde{u}|^p{\mathrm{d}}x,\left( \int_{\Omega}\Phi^p|\nabla \tilde{u}|^p{\mathrm{d}}x\right)^\frac{2}{p}\right\}\,.
\end{equation*} 
Combining case 1 and case 2, we obtain \eqref{rem...}. Hence the proof.
\end{proof}
 Now we prove  quantitative stability of \eqref{HS} using Lemma \ref{Poincare}, Proposition \ref{prop-radi-dec}, and Lemma \ref{reminder-lemma}. First, we prove Theorem \ref{Hardy-Thm-p}.

\begin{proof}[Proof of Theorem \ref{Hardy-Thm-p}]
    Let $\Omega=B_R(x_0),\,R>0.$
 Due to the homogeneity-free formulation of \eqref{Eqn-Hardy-p}, it is enough to prove 
\begin{align} \label{Eq:main_req}
    \delta_p(u)\gtrsim D_p(u)^{\max\{2,p\}}
\end{align}
for all $u\in \Dp$ with $\|\nabla u\|_{L^p(\Omega)}=1.$ Moreover, $C_c^1(\Omega)$ is dense in $\Dp.$ Thus, it is sufficient to prove \eqref{Eq:main_req} for all $u\in C_c^1(\Omega)$ satisfying $\|\nabla u\|_{L^p(\Omega)}=1.$ 
Choose such $u$ arbitrary. Note that the left-hand side of \eqref{Eq:main_req} is
$$\delta_p(u)=\int_\Omega|\nabla u|^p\dd x- \mathbb{H}_p(\Omega)^p\int_\Omega|g||u|^p \ {\mathrm{d}}x\,.
$$ Thus, applying Lemma \ref{reminder-lemma}, we obtain the estimate 
  \begin{align}\label{deficit-estimate1}   
\delta_p(u)\ge C 
\begin{cases}
\int_\Omega \left|\nabla \left(\frac{u}{\Phi}\right) \right|^p\Phi^p\,{\mathrm{d}}x
& \text{if } p\ge 2, \\[1ex]
\min\left\{\int_\Omega \left|\nabla \left(\frac{u}{\Phi}\right) \right|^p\Phi^p\,{\mathrm{d}}x,\,\left[ \int_\Omega \left|\nabla \left(\frac{u}{\Phi}\right) \right|^p\Phi^p\,{\mathrm{d}}x\right]^\frac{2}{p}\right\}\
& \text{if}\; p\in (1,2)
\end{cases}
\end{align}
for some $C=C(p)>0.$ 
 Next, we estimate the integral on the right-hand side of \eqref{deficit-estimate1} using a weighted Poincar\'e inequality (see Lemma \ref{Poincare}). Thanks to Proposition \ref{prop-radi-dec}, which proves that $\Phi$ is radially nonincreasing.  Thus, we apply the weighted Poincar\'e inequality in Lemma \ref{Poincare} with weight $\Phi^p$ to obtain
\begin{align*}
     \int_\Omega \left|\nabla u\right|^p\Phi^p\,{\mathrm{d}}x& \ge C_1\int_\Omega|u-(u)_{\Phi^p}|^p \,\Phi^p\, {\mathrm{d}}x \ge C_1 \inf_{c\in \mathbb{R}}\int_\Omega|u-c|^p \,\Phi^p\, {\mathrm{d}}x,
\end{align*}
for some $C_1=C_1(N,p,R,\Phi)>0$. Now, replace $u$ with $\frac{u}{\Phi}$ to obtain 
\begin{align}
     \int_\Omega \left|\nabla \left(\frac{u}{\Phi}\right)\right|^p\Phi^p\,{\mathrm{d}}x&\ge C_1\inf_{c\in \mathbb{R}}\int_\Omega\left|\frac{u}{\Phi}-c\right|^p\,\Phi^p\,{\mathrm{d}}x \nonumber\\&=C_1\inf_{c\in \mathbb{R}}\int_\Omega\left|u-c\,\Phi\right|^p{\mathrm{d}}x=C_1D_p(u)^p\label{Eqn1-thm-p-Hardy}.
\end{align}
If $p\ge 2$ then the stability estimate \eqref{Eq:main_req} follows by combining \eqref{deficit-estimate1} and
 \eqref{Eqn1-thm-p-Hardy}. Now assume that $p\in (1,2).$ In this case, combining \eqref{deficit-estimate1} and
 \eqref{Eqn1-thm-p-Hardy}, we obtain 
 \begin{align}
    \delta_p(u)\ge C_2\,\min\{D_p(u)^p, D_p(u)^2\}\label{Eq8<2}.
 \end{align}
 for some $C_2=C_2(N,p,R,\Phi)>0$. Now, by the definition of $D_p(u),$ we obtain
$$D_p(u)\le \frac{\|u\|_{L^p(\Omega)}}{\|\nabla u\|_{L^p(\Omega)}}\le C_3,$$ for some $C_3=C_3(N,p,R)>0,$ where the last inequality follows by using the Poincar\'e inequality. That is, $D_p(u)/C_3\le 1.$ Using this in \eqref{Eq8<2}, we obtain  
 \begin{align*}
    \delta_p(u)\ge C_4\, D_p(u)^2\quad\text{for}\;p\in (1,2),
 \end{align*}
for some constant $C_4=C_4(N,p,R,\Phi)>0.$ That is \eqref{Eq:main_req} holds, completing the proof. 
\end{proof}

Next, we prove the gradient stability of \eqref{HS}, as stated in Theorem \ref{thm-gard>2}.
\begin{proof}[Proof of Theorem \ref{thm-gard>2}] 

 \textbf{ Stability : Proof of \eqref{Eqn-Hardy-p-infty}.} It is enough to prove the stability inequality \eqref{Eqn-Hardy-p-infty} for all $u\in \Dp$ satisfying $\|\nabla u\|_{L^p(\Omega)}=1.$ Further, since $ C_c^\infty(\Omega)$ is dense in $\Dp,$ it is sufficient to prove \eqref{Eqn-Hardy-p-infty}
for all $u\in  C_c^\infty(\Omega)$ satisfying  $\|\nabla u\|_{L^p(\Omega)}=1.$ Choose such $u$ arbitrary. 
Applying Lemma \ref{reminder-lemma}, we obtain 
  \begin{align}\label{deficit-estimate2}   
\delta_p(u)\ge C 
\begin{cases}
[I(u,\Phi)]^p
& \text{if } p\ge 2, \\[1ex]
\min\left\{[I(u,\Phi)]^p,\,\left[ I(u,\Phi)\right]^2\right\}\
& \text{if}\; p\in (1,2)
\end{cases}
\end{align} 
for some $C=C(p)>0,$ where $$I(u,\Phi)=\left(\int_\Omega \left|\nabla \left(\frac{u}{\Phi}\right) \right|^p\Phi^p\,{\mathrm{d}}x\right)^{1/p}\,.$$
Next, we prove the estimate $I(u,\Phi)\gtrsim d_p(u)$ using Proposition \ref{prop-radi-dec} and the weighted Poincar\'e inequality in Lemma \ref{Poincare}. We adapt the idea of the proof of \cite[Theorem 3.1]{Amir2} to establish this estimate. 
Let $c_0=c_0(u)\in \R$ be such that 
\begin{equation}
\int_\Omega u\Phi^{p-1}\dd x-c_0\int_\Omega\Phi^p\dd x=0.    
\end{equation}
Define $$f:=\frac{u}{\Phi}-c_0.$$ For this choice of $c_0$, we have $(f)_{\Phi^p}=\int_\Omega f\Phi^p\dd x=0.$ Moreover, $f\in W^{1,p}(\Omega)$ since $u\in C_c^\infty(\Omega),$ $\Phi\in C_{\mathrm{loc}}^1(\Omega),$ and $\Phi>0$ in $\Omega$ (see Lemma \ref{cuesta-lem}). Thus, using the weighted Poincar\'e inequality in Lemma \ref{Poincare} with radially nonincreasing weight $\Phi^p$ (Proposition \ref{prop-radi-dec}), we obtain 
\begin{equation}\label{poincare-estimate-infty}
  \int_\Omega|\nabla f|^p\Phi^p\dd x\ge C_1 \int_\Omega|f|^p\Phi^p\dd x\,,
\end{equation}
for some $C_1=C_1(N,p,R,\Phi)>0.$ 
Observe that
\begin{align}\label{estimate-distance2}
    d_p(u)=\inf\left\{
\frac{\|\nabla (u-c\Phi)\|_{L^p(\Omega)}}
{\|\nabla u\|_{L^p(\Omega)}}:
c\in\R
\right\}&\le \|\nabla (u-c_0\Phi)\|_{L^p(\Omega)}=\|\nabla (\Phi f)\|_{L^p(\Omega)}\nonumber\\&\le\|\Phi\nabla f\|_{L^p(\Omega)}+\|f \nabla \Phi\|_{L^p(\Omega)},
\end{align}
where the last inequality follows by the triangular inequality. Note that 
\begin{equation}\label{notation}
    \|\Phi\nabla f\|_{L^p(\Omega)}=I(u,\Phi).
\end{equation}
So, we only need to estimate $\|f \nabla \Phi\|_{L^p(\Omega)}.$
 Notice that $f$ and $\nabla f$ are bounded in compact set $K$ containing the support of $u$ and $f=-c_0$ on $\Omega\setminus K$, consequently, $f$ and $\nabla f$ are bounded in $\Omega.$  Moreover, product rule gives 
$$\nabla (\Phi|f|^p)=|f|^p\nabla \Phi+p\Phi|f|^{p-2}f\nabla f.$$ Therefore, $\nabla (\Phi|f|^p)\in L^p(\Omega)$ and $\Phi|f|^p\in W_0^{1,p}(\Omega).$ Thus, we can choose $v=\Phi|f|^p$ as a test function in the weak formulation \eqref{Euler-Lag-p} to get
\begin{equation*}
    \int_\Omega |f|^p|\nabla \Phi|^{p}\dd x+p\int_\Omega\Phi|f|^{p-2}f|\nabla \Phi|^{p-2}\nabla\Phi\cdot\nabla f\dd x=\mathbb{H}_{p}(\Omega)^p\,\int_\Omega |g|\Phi^{p}|f|^p\,\dd x.
\end{equation*}
Since $g\in L^\infty(\Omega),$ the above inequality implies that 
\begin{equation*}
    \int_\Omega |f|^p|\nabla \Phi|^{p}\dd x\le \mathbb{H}_{p}(\Omega)^p\|g\|_{L^\infty(\Omega)}\,\int_\Omega \Phi^{p}|f|^p\,\dd x-p\int_\Omega\Phi|f|^{p-2}f|\nabla \Phi|^{p-2}\nabla\Phi\cdot\nabla f\dd x.
\end{equation*}
The H\"older inequality yields
$$\int_\Omega\Phi|f|^{p-2}f|\nabla \Phi|^{p-2}\nabla\Phi\cdot\nabla f\dd x\le \left(\int_\Omega|f\nabla \Phi|^p\dd x\right)^{(p-1)/p} \left(\int_\Omega|\Phi\nabla f|^p\dd x\right)^{1/p}.$$ Consequently, using the Poincar\'e estimate \eqref{poincare-estimate-infty}, we obtain  
$$ \int_\Omega |f\nabla \Phi|^{p}\dd x\le C_2\,\int_\Omega |\nabla f|^p\Phi^p\,\dd x+p\left(\int_\Omega|f\nabla \Phi|^p\dd x\right)^{(p-1)/p}\left(\int_\Omega|\Phi\nabla f|^p\dd x\right)^{1/p},$$ where $C_2=C_1\,\mathbb{H}_{p}(\Omega)^p\|g\|_{L^\infty(\Omega)}.$ Defining $$J(u,\Phi):=\|f\nabla\Phi\|_{L^p(\Omega)} .$$ 
and using $ I(u,\Phi)=\|\Phi\nabla f\|_{L^p(\Omega)}$ (see \eqref{notation}),  the above inequality can be rewritten as 
\begin{equation}\label{E1:i-inf}
   [ J(u,\Phi)]^p\le C_2\,[I(u,\Phi)]^p+p\,[J(u,\Phi)]^{p-1}\,I(u,\Phi)\,.
\end{equation} 
Now we show that $J(u,\Phi)\le cI(u,\Phi)$ for some $c>0$ independent of $u.$ If $I(u,\Phi)=0$ then by \eqref{E1:i-inf} we get $J(u,\Phi)=0,$ so assume that $I(u,\Phi)>0.$ Dividing both sides of \eqref{E1:i-inf} by $[I(u,\Phi)]^p$ and arranging we obtain 
 \begin{equation}\label{E2:i-inf}
    \left[\frac{J(u,\Phi)}{I(u,\Phi)}\right]^{p-1}\,\left[\frac{J(u,\Phi)}{I(u,\Phi)}-p\right]\le C_2.
\end{equation} We know that
\begin{equation*}
  L: t\mapsto t^{p-1}(t-p)-C_2
\end{equation*}
is an increasing continuous function on $[p,\infty) ,$  $L(p)<0,$ and $L(t)>0$ for large $t>p.$ Thus, by the Intermediate value theorem,  there exists a zero $\gamma$ in $(p,\infty)$, and it is unique since $L$ is strictly increasing. If $J(u,\Phi)/I(u,\Phi)\le p,$ then $J(u,\Phi)/I(u,\Phi)\le \gamma$. On the other hand, $L[J(u,\Phi)/I(u,\Phi)]\le 0=L(\gamma)$ (see \eqref{E2:i-inf}) which implies that  $J(u,\Phi)/I(u,\Phi)\le \gamma,$ since $L$ is strictly increasing. Thus, in both cases, $J(u,\Phi)\le \gamma I(u,\Phi),$ i.e., 
\begin{equation}\label{E3:i-inf}
 \|f\nabla \Phi\|_{L^p(\Omega)} \le \gamma\|\Phi\nabla f \|_{L^p(\Omega)}\,.
\end{equation} 
Here $\gamma$ depends only on $p$ and the constant $C_2.$
Combining \eqref{estimate-distance2} and \eqref{E3:i-inf}, we obtain 
\begin{equation*}
    d_p(u)\le \|\Phi\nabla f\|_{L^p(\Omega)}+\gamma\|\Phi\nabla f\|_{L^p(\Omega)}=  (1+\gamma) I(u,\Phi).
\end{equation*}
Thus, using the above estimate in \eqref{deficit-estimate2}, we obtain 
\[\delta_p(u)\ge \frac{C}{(1+\gamma)^p} 
\begin{cases}
d_p(u)^p
& \text{if } p\ge 2, \\[1ex]
\min\left\{d_p(u)^p, (1+\gamma)^{p-2}\,d_p(u)^2\right\}\ge (1+\gamma)^{p-2} d_p(u)^2
& \text{if}\; p\in (1,2).
\end{cases}\] 
where we used the fact that $(1+\gamma)^{p-2}\le 1$ for $p\in (1,2)$ and $d_p(u)\le 1.$ That is, the stability estimate \eqref{Eqn-Hardy-p-infty} holds. Hence the proof.

\smallskip

\textbf{Optimality of the exponent $\alpha=\max\{2,p\}$:} We follow the proof of Theorem 3.1 in \cite{Amir2} to prove the optimality. Suppose that for $\beta\in (0,\alpha),$ there exists a constant $C>0$ such that
\begin{equation}\label{E1-opt}
    \delta_p(u)\ge Cd_p(u)^\beta\quad\forall\,u\in \Dp.
\end{equation} 
We construct perturbations of $\Phi$ that contradict \eqref{E1-opt}. 
Since $\|\nabla \Phi\|_{L^p(\Omega)}>0$ and $\Phi\in C_{\mathrm{loc}}^1(\Omega)$ (see Lemma \ref{cuesta-lem}), there exists $x_0\in \Omega$ such that $\nabla\Phi(x_0)\ne 0.$ Now continuity implies that there exists a sufficiently small $r>0$ satisfying 
\begin{equation}\label{E2-opt}
    0<m\le|\nabla\Phi|\le M\quad\text{on}\quad B_r(x_0)\qquad\text{and}\quad M^p|B_1(0)|r^N\le \frac{1}{2}\|\nabla \Phi\|_{L^p(\Omega)}^p.
\end{equation}
for some $M\ge m>0.$ Set $B:=B_r(x_0).$ Use \eqref{E2-opt} to obtain
\begin{equation}\label{E3-opt*}
    \|\nabla \Phi\|_{L^p(\Omega\setminus B)}^p=\|\nabla \Phi\|_{L^p(\Omega)}^p-\|\nabla \Phi\|_{L^p(B)}^p\ge\|\nabla \Phi\|_{L^p(\Omega)}^p-M^p|B_1(0)|r^N\ge \frac{1}{2}\|\nabla \Phi\|_{L^p(\Omega)}^p.
\end{equation}
Choose $h\in C_c^\infty(\Omega)$ with the support contained in $B.$ Then $\nabla h=0$ on $\Omega\setminus B,$ and hence
\begin{align}\label{E3-opt}
    \|\nabla (\Phi+h-a\Phi)\|_{L^p(\Omega)}^p\ge \int_{\Omega\setminus B}|\nabla (\Phi-a\Phi)|^p\dd x&\ge |1-a|^p\|\nabla \Phi\|^p_{L^p(\Omega\setminus B)}\nonumber\\&\ge\frac{1}{2}|1-a|^p\|\nabla \Phi\|^p_{L^p(\Omega)},\quad \forall a\in\R
\end{align}
where the last inequality follows by \ref{E3-opt*}. Since $\nabla h=\nabla(\Phi +h-a\Phi)-(1-a)\nabla \Phi,$ using the triangular inequality and \eqref{E3-opt}, we obtain
\begin{align*}
    \|\nabla h\|_{L^p(\Omega)}&\le \|\nabla (\Phi+h-a\Phi)\|_{L^p(\Omega)}+|1-a|\|\nabla \Phi\|_{L^p(\Omega)}\\&\le(1+2^{1/p}) \|\nabla (\Phi+h)-a\nabla\Phi\|_{L^p(\Omega)}\qquad \forall a\in \R.
\end{align*}
Dividing both sides of the above inequality by $\|\nabla(\Phi+h)\|_{L^p(\Omega)}$ and then taking the infimum over $a\in \R$ yields 
\begin{equation}\label{E4-opt}
\frac{\|\nabla h\|_{L^p(\Omega)}}{\|\nabla(\Phi+h)\|_{L^p(\Omega)}}\le (1+2^{1/p})\,d_p(\Phi+h).
\end{equation}
 Next, we estimate the deficit $\delta_p(\Phi+h).$ Define the convexity reminder
 $$\Gamma_p(s,t):=|s|^p-|s-t|^p-p|s-t|^{p-2}(s-t)t\ge 0.$$
  Then $$\Gamma_p(\Phi+h,h)=|\Phi+h|^p-\Phi^p-p\Phi^{p-1}h.$$ Taking $h$ as a test function in the weak formulation \eqref{Euler-Lag-p}, we obtain 
  \begin{equation}\label{E4-otp**}
  \int_B|\nabla\Phi|^{p-2}\nabla\Phi\cdot\nabla h\dd x=\mathbb{H}_p(\Omega)^p\int_B|g|\Phi^{p-1}h\dd x.    
  \end{equation}
   Thus, using $\Phi\in \mathcal{M}_p,$ $|\Phi+ h|^p-|\Phi|^p=0$ on $B^c,$ and  \eqref{E4-otp**} we obtain
  \begin{align} \|\nabla(\Phi+h)\|_{L^p(\Omega)}^p\,&\delta_p(\Phi+h)=  \int_\Omega  |\nabla\Phi+\nabla h|^p \dd x-\mathbb{H}_p(\Omega)^p\int_\Omega|g||\Phi+h|^p\dd x\nonumber\\&=  \int_B\left(|\nabla\Phi+\nabla h|^p-|\nabla\Phi|^p\right)\dd x-\mathbb{H}_p(\Omega)^p\int_B|g|(|\Phi+h|^p-|\Phi|^p)\dd x\nonumber\\&=  \int_B\Gamma_p(\nabla\Phi+\nabla h,\nabla h)\dd x-\mathbb{H}_p(\Omega)^p\int_B|g|\,\Gamma_p(\Phi+h,h)\dd x \nonumber\\&\le  \int_B\Gamma_p(\nabla\Phi+\nabla h,\nabla h)\dd x\,,\quad\text{since}\;\;\Gamma_p\ge 0.\label{E5-opt}
  \end{align}

Case 1: $p\in (1,2].$ Fix a nonzero function $\Psi\in C_c^\infty(B)$ and consider $\Phi_t:=\Phi+t\Psi,\;t\in \R.$ We estimate $d_p(\Phi_t)$ and $\delta_p(\Phi_t)$ by using \eqref{E4-opt} and \eqref{E5-opt}. Take $h=t\Psi$ in \eqref{E4-opt} to get 
\begin{equation}\label{E4-opt-modifty1}
|t|  \frac{\|\nabla \Psi\|_{L^p(\Omega)}}{\|\nabla\Phi_t\|_{L^p(\Omega)}}\le (1+2^{1/p})\, d_p(\Phi_t)\,.
\end{equation}
Observe that 
\begin{equation}\label{Eq6-opt}
   \frac{1}{2} \|\nabla\Phi\|_{L^p(\Omega)}\le \|\nabla\Phi_t\|_{L^p(\Omega)}\le \frac{3}{2}\|\nabla\Phi\|_{L^p(\Omega)}\;\;\text{for sufficiently small}\;|t|. 
\end{equation}
Combining \eqref{E4-opt-modifty1} and \eqref{Eq6-opt} we get 
  $$d_p(\Phi_t)\gtrsim |t|\quad \text{for sufficiently small}\;|t|.$$ 
 Since $t^{p/2}$ is concave on $[0,\infty),$ 
 we also have 
\begin{align*}
|\nabla\Phi+t\nabla\Psi|^p&=\left[|\nabla \Phi|^2+2t\nabla\Phi\cdot\nabla \Psi+t^2|\nabla\Psi|^2\right]^{p/2}\\&\le |\nabla \Phi|^p+\frac{p}{2}|\nabla \Phi|^{p-2}\left( 2t\nabla\Phi\cdot\nabla \Psi+t^2|\nabla\Psi|^2\right)
\end{align*}
which implies that 
\begin{align*}
    \Gamma_p(\nabla \Phi+t\nabla\Psi,t\nabla\Psi)&=|\nabla\Phi+t\nabla\Psi|^p-|\nabla\Phi|^p-pt|\nabla\Phi|^{p-2}\nabla\Phi\cdot\nabla\Psi\\&\le \frac{p}{2}t^2|\nabla \Phi|^{p-2}|\nabla\Psi|^2\le \frac{p}{2}t^2m^{p-2}|\nabla\Psi|^2\quad \text{on}\; B
\end{align*}
where the last inequality follows using  \eqref{E2-opt} ($|\nabla\Phi|\ge m$ on $B$) and $p-2\le 0.$ 
 Use this estimate in \eqref{E5-opt} (with  $h=t\Psi$) to obtain  
  $$\delta_p(\Phi_t)\le \frac{1}{\|\nabla \Phi_t\|_{L^p(\Omega)}^p} \frac{p}{2}m^{p-2}t^2\int_B|\nabla\Psi|^2\dd x\le p2^{p-1}m^{p-2}\frac{\|\nabla \Psi\|_{L^2(\Omega)}^2}{\|\nabla \Phi\|_{L^p(\Omega)}^p}\,t^2\,,$$  where the last inequality follows from  \eqref{Eq6-opt}. That is $$ \delta_p(\Phi_t)\lesssim t^2\quad \text{for sufficiently small}\;|t|.$$ Therefore, applying \eqref{E1-opt} to $\Phi_t$, we obtain $t^2\gtrsim\delta_p(\Phi_t)\ge C d_p(\Phi_t)^\beta\gtrsim |t|^\beta$ for sufficiently small $|t|$. Since $\beta<2$, this gives a contradiction.

\smallskip

  Case 2: $p>2.$ Fix a nonzero $\Psi\in C_c^\infty(B_1(0))$ and let $0<\epsilon<r$.
Define
$$h_\epsilon(x):=\epsilon\Psi\left(\frac{x-x_0}{\epsilon}\right)\quad\text{and}\quad \Phi_\epsilon:=\Phi+h_\epsilon\,.$$ Then we have   $$\text{supp}\,h_\epsilon\,\subset B_\epsilon(x_0)\subset B\quad\text{and}\quad\nabla h_\epsilon(x)=\nabla\Psi\left(\frac{x-x_0}{\epsilon}\right).$$ Therefore, making the change of variables formula $(y=(x-x_0)/\epsilon)$, we can obtain $$\|\nabla h_\epsilon\|_{L^p(\Omega)}=\epsilon^{N/p}\|\nabla\Psi\|_{L^p(B_1(0))}\,.
$$ We estimate $\delta_p(\Phi_t)$ and $d_p(\Phi_t)$ in the same way as in case $p\in (1,2].$ 
Substitute the above equality in \eqref{E4-opt} to get \begin{equation*}
  \epsilon^{N/p}\frac{\|\nabla \Psi\|_{L^p(\Omega)}}{\|\nabla\Phi_\epsilon\|_{L^p(\Omega)}}\le (1+2^{1/p})\,d_p(\Phi_\epsilon)\,.
\end{equation*} Moreover, it can be verified that 
\begin{equation}\label{Eq10-opt}
   \frac{1}{2} \|\nabla\Phi\|_{L^p(\Omega)}\le \|\nabla\Phi_\epsilon\|_{L^p(\Omega)}\le \frac{3}{2}\|\nabla\Phi\|_{L^p(\Omega)}\;\;\text{for sufficiently small}\;\epsilon. 
\end{equation}
Therefore, $$d_p(\Phi_\epsilon)\gtrsim \epsilon^{N/p}\quad \text{for sufficiently small}\;\epsilon.$$
Notice that $|\nabla h_\epsilon|\le \|\nabla\Psi\|_{L^\infty(B_1(0))}$ and $|\nabla\Phi|<M$ on $B.$ Thus, 
\begin{align*}
\Gamma_p(\nabla\Phi+\nabla h_\epsilon,\nabla h_\epsilon)&=|\nabla\Phi+\nabla h_\epsilon|^p-|\nabla h_\epsilon|^p-p|\nabla\Phi|^{p-2}\nabla\Phi\cdot\nabla h_\epsilon\nonumber\\&\le \left(|\nabla\Phi|+|\nabla h_\epsilon|\right)^p+p|\nabla\Phi|^{p-1}|\nabla h_\epsilon|\nonumber\\&\le  \left(M+\|\nabla\Psi\|_{L^\infty(B_1(0))}\right)^p+pM^{p-1}\|\nabla\Psi\|_{L^\infty(B_1(0))}=C^*\,,  
\end{align*}
where $C^*>0$ is independent of $\epsilon.$
Which gives $$\int_B\Gamma_p(\nabla\Phi+\nabla h_\epsilon,\nabla h_\epsilon)\dd x=\int_{B_\epsilon}\Gamma_p(\nabla\Phi+\nabla h_\epsilon,\nabla h_\epsilon)\dd\le C^*
|B_1(0)|\epsilon^N
.$$
Combining this with \eqref{E5-opt} and \eqref{Eq10-opt} we obtain $$\delta_p(\Phi_\epsilon)\lesssim\epsilon^N \quad \text{for sufficiently small}\;\epsilon.$$  Thus, applying \eqref{E1-opt} to $\Phi_\epsilon$, we obtain $\epsilon^N\gtrsim \epsilon^{\beta N/p}$. Since $p>\beta$, this leads to a contradiction as $\epsilon\to0$. Combining Cases 1 and 2, we conclude that the exponent $\alpha=\max\{2,p\}$ in \eqref{Eqn-Hardy-p-infty} is optimal. This completes the proof.
 \end{proof}
\begin{rem}\label{last-rem} As seen in the proof above, establishing the sharpness of the exponent $\alpha=\max\{2,p\}$ in \eqref{Eqn-Hardy-p-infty} requires only local $C^1$-regularity of the extremal $\Phi$ (see \eqref{E2-opt}). Therefore, by using the same argument, we can show that the exponent $\alpha=2$ in the gradient stability estimate \eqref{Eqn-sta-2} (of Theorem \ref{Thm-p=2}) is also sharp under the same regularity assumption. This is discussed in Remark \ref{rem---sharp=2}.

\end{rem}
\textbf{Question:} 
 For $p\neq 2$, we prove quantitative gradient stability for \eqref{HS} only when the Hardy potentials $g$ is radial and belongs to $L^\infty(\Omega)$ (see Theorem \ref{thm-gard>2}). Can we extend this result to unbounded Hardy potentials $g$, particularly those considered in Theorem \ref{Hardy-Thm-p} ?

  \begin{center}
	{\bf Acknowledgement}
\end{center}
S Roy, acknowledges the support of the National Board for Higher Mathematics (NBHM) through the NBHM Postdoctoral Fellowship
 (No.0202/3/2026/R\&D-II/2870). The author thanks Prof. Ujjal Das for valuable discussion on the article.

\bibliographystyle{abbrvurl}
\bibliography{Reference}
\end{document}